\documentstyle{amsppt}
\magnification=1200
\NoBlackBoxes

\topmatter
\title On the embedding of 2-concave
Orlicz spaces into $L^1$
\endtitle
\author
Carsten Sch\"utt
\endauthor
\address
Oklahoma State University, Department of Mathematics,
Stillwater, OK 74078 
\endaddress
\address
Christian Albrechts Universit\"at, Mathematisches Seminar,
24118 Kiel, Germany
\endaddress
\abstract
In [K--S 1] it was shown that

$$ 
\underset {\pi} \to {\text{Ave}}
(\sum_{i=1}^{n}|x_i a_{\pi(i)}|^2)^{\frac {1}{2}}
$$

is equivalent to an Orlicz norm whose Orlicz function
is 2-concave. Here we give a formula for the sequence
$a_1, a_2,....,a_n$ so that the above expression is
equivalent to a given Orlicz norm.
\endabstract
\thanks
Supported by NSF-grant DMS-9301506
\endthanks
\endtopmatter
\pagebreak

\document

A convex function $M\: \Bbb R \to \Bbb R$ with $M(t)=M(-t)$,
$M(0)=0$, and $M(t)>0$ if $t \neq 0$ is called an Orlicz
function. M is said to be 2-concave if $M(\sqrt t)$ is a 
concave function on $[0, \infty)$ and strictly 2-concave
if $M(\sqrt t)$ is strictly concave. M is 2-convex if $M(\sqrt{t})$
is convex and strictly 2-convex if $M(\sqrt t)$ is
strictly convex. If $M'$ is invertible
on $(0, \infty)$ then the dual function is given by

$$
M^{*}(t)=\int_{0}^{t} {M'}^{-1}(s)ds
$$

We define the Orlicz norm of a sequence $\{x_{i}\}_{i=1}^{\infty}$
by

$$
\parallel x \parallel_{M}=\text{sup} \{ \sum_{i=1}^{\infty}x_iy_i
 |\sum_{i=1}^{\infty} M^{*}(y_i) \leq 1 \}
$$

In [K--S 1, K--S 2] we have used a different expression for the
definition of the Orlicz norm: x has norm equal to 1 if and only
if $\sum_{i=1}^{\infty} M(x_i)=1$. But it turns out that the above
definition gives slightly better estimates. \par
Bretagnolle and Dacunha-Castelle [B--D] showed that an Orlicz space
$l^M$ is isomorphic to a subspace of $L^1$ if and only if M is
equivalent to a 2-concave Orlicz function. As a corollary we get the
same result here. In [K--S 1] a variant of
the following
result was obtained.
\vskip 1cm

\proclaim{\smc Theorem 1}
Let $a_1 \geq a_2 \geq ... \geq a_n > 0$ and let M be the Orlicz
with

$$
{M^{*}}^{-1}(\frac{l}{n})
= \left\{ (\frac{1}{n} \sum_{i=1}^{l} a_{i})^{2} +
\frac{l}{n}
(\frac{1}{n} \sum_{i=l+1}^{n}|a_{i}|^{2})\right\}^{\frac{1}{2}}
\tag 1
$$

for all $l=1,2,...,n$ and  such that ${M^{*}}^{-1}$ is an affine
function between the given values. Then we have for all $x \in \Bbb R^n$

$$ 
\frac{1}{2\sqrt{5}}\frac{(n-1)^2}{n^2+(n-1)^2}
 \parallel x \parallel_{M} \leq
\underset {\pi} \to {\text{Ave}}
(\sum_{i=1}^{n}|x_i a_{\pi(i)}|^2)^{\frac {1}{2}} \leq
\frac{2\sqrt{2}}{c_{n}}\parallel x \parallel_{M}
\tag 2
$$

where $c_n=1-\frac{1}{2!}+\frac{1}{3!}-...+(-1)^{n+1}\frac{1}{n!}$.
\endproclaim
\vskip 1cm

We present here those arguments of the proof of Theorem 1 that are
different from the arguments in [K--S 1, K--S 2].  \par
There is always an Orlicz function M satisfying the assumptions of
Theorem 1. In order to verify this we show that

$$
\left \{ (\int_{0}^{x} f(t)dt)^{2}+x \int_{x}^{1}|f(t)|^{2}dt
\right \}^{\frac{1}{2}}
$$

is a concave function of x. Moreover, we may assume that f is
differentiable. For the second derivative of the above expression
we get

$$
\frac{f'(x)(\int_{0}^{x} f(t)dt-xf(x))}
{\left \{ (\int_{0}^{x} f(t)dt)^{2}+x \int_{x}^{1}|f(t)|^{2}dt 
\right \}^{\frac{1}{2}}}
-\frac{1}{4} \frac{(2f(x)\int_{0}^{x} f(t)dt +\int_{x}^{1}|f(t)|^{2}dt
-xf(x))^{2}}{\left \{ (\int_{0}^{x} f(t)dt)^{2}+x \int_{x}^{1}|f(t)|^{2}dt  
\right \}^{\frac{3}{2}}}
$$

The first summand is nonpositive since f is decreasing. \par

It follows from Theorem 1 that an Orlicz function M has to
be equivalent to a 2-concave Orlicz function if $l^{M}$ is
isomorphic to a subspace of $L^{1}$ [K--S 1, K--S 2]. We compute
here how we have to choose the sequence $a_1, a_2,..., a_n$ so
that we get (1) for a given 2-concave Orlicz function M. From this
it also follows that $l^{M}$ is isomorphic to a subspace of $L^{1}$
if M is 2-concave.
\vskip 1cm

\proclaim{\smc Theorem 2}
Let M be a strictly convex, twice differentiable Orlicz function that
is strictly 2-concave. Assume that
${M^{*}}(1)=1$ and let

$$
a_{l}= -\frac{n}{2} \int_{\frac{l-1}{n}}^{\frac{l}{n}}
\int_{t}^{1} \frac{ (({M^{*}}^{-1})^2)''(s) }{\sqrt{ ({M^{*}}^{-1})^{2}(s)-s
(({M^{*}}^{-1})^2)'(s)}} ds +1-\sqrt{1-(({M^{*}}^{-1})^2)'(1)} dt
\tag 3
$$

for $l=1,2,...,n$. Then we have for all
$x \in \Bbb R^n$

$$
\frac{1}{c} \parallel x \parallel_{M} \leq
\underset {\pi} \to {\text{Ave}}
(\sum_{i=1}^{n}|x_i a_{\pi(i)}|^2)^{\frac {1}{2}} \leq
c \parallel x \parallel_{M}
$$

where c is a constant that does not depend on n and M.
\endproclaim
\vskip 1cm

Since $\Bbb R^n$ with the norm

$$
\parallel x \parallel =\underset {\pi} \to {\text{Ave}}
(\sum_{i=1}^{n}|x_i a_{\pi(i)}|^2)^{\frac {1}{2}}
$$

is isometric to a subspace of $L^1$, we get the following corollary.
\vskip 1cm

\proclaim{\smc Corollary 3}
Let M be a 2-concave Orlicz function. Then $l^M$ is isomorphic to
a subspace of $L^1$.
\endproclaim
\vskip 1cm

\proclaim{\smc Lemma 4}[K--S 1]
For all $n \in \Bbb N$ and all $n \times n$ matrices A with nonnegative
entries we have

$$
c_n \frac{1}{n} \sum_{k=1}^{n} s(k) \leq
\underset {\pi} \to {\text{Ave}} \underset {1 \leq i \leq n} \to {\text{max}}
|a(i,\pi(i))| \leq \frac{1}{n} \sum_{k=1}^{n} s(k)
$$

where $c_n=1-\frac{1}{2!}+\frac{1}{3!}-...+(-1)^{n+1}\frac{1}{n!}$ and $s(k),
k=1,2,...,n^2$ is the nonincreasing rearrangement of the numbers $a(i,j),
i,j=1,2,...,n$.

\endproclaim
\vskip 1cm

\proclaim{\smc Lemma 5}[K--S 2]
For all $n \in \Bbb N$ and all
nonnegative numbers $a(i,j,k), i,j,k=1,2,...,n$ we
have

$$
\frac{(n-1)^2}{n^2+(n-1)^2} \frac{1}{n^2} \sum_{k=1}^{n^2} s(k) \leq
\underset {\pi,\sigma} \to {\text{Ave}} \underset{1 \leq i \leq n} \to
{\text{max}} |a(i,\pi(i),\sigma(i))| \leq
 \frac{1}{n^2} \sum_{k=1}^{n^2} s(k)
$$

where the average is taken over all permutations $\pi, \sigma$ of
$\{1,2,...,n\}$ and $\{s(k)\}_{k=1}^{n^3}$ is the nonincreasing
rearrangement of the numbers $a(i,j,k), i,j,k=1,2,...,n$.
\endproclaim
\vskip 1cm

\proclaim{\smc Lemma 6}[K--S 1]
Let $b_1 \geq b_2 \geq ... \geq b_s > 0$, $n \leq s$, and 

$$
\parallel x \parallel_{b} =
\underset {\sum k_{j}=s} \to {\text{max}}
\sum_{i=1}^{n}(\sum_{j=1}^{k_i} b_j)|x_{i}|
$$

Then we have for all Orlicz functions M with
$M^{*}(\sum_{j=1}^{l} b_j)=\frac{l}{s}, l=1,...,s$ and all
$x \in \Bbb R^n$

$$
\parallel x \parallel_{b} \leq
\parallel x \parallel_{M} \leq
2\parallel x \parallel_{b}  
$$

\endproclaim
\vskip 1cm

The proof of the right hand inequality of Lemma 6 is the same as
in [K--S 2]. The left hand inequality follows from the definition
of the norm.

\vskip 1cm
\demo{Proof of Theorem 1}
We choose the sequence $b_j, j=1,2,...,n$ with

$$
\sum_{j=1}^{k} b_j= \sqrt {nk}  \qquad k=1,2,...,n 
$$

Then we get by Lemmata 4 and 6

$$
\frac{c_n}{2} \underset {\pi} \to {\text{Ave}}
(\sum_{i=1}^{n}|x_i a_{\pi(i)}|^2)^{\frac {1}{2}} \leq
\underset {\pi,\sigma} \to {\text{Ave}} \underset{1 \leq i \leq n} \to
{\text{max}} |x_{i}a_{\pi(i)}b_{\sigma(i)}| \leq
\underset {\pi} \to {\text{Ave}}
(\sum_{i=1}^{n}|x_i a_{\pi(i)}|^2)^{\frac {1}{2}}
$$

And by Lemma 5 we get

$$
\frac{(n-1)^2}{n^2+(n-1)^2}  \frac{1}{n^2}\sum_{k=1}^{n^2} s(k) \leq
(\sum_{i=1}^{n}|x_i a_{\pi(i)}|^2)^{\frac {1}{2}} \leq
\frac{2}{c_n} \frac{1}{n^2}\sum_{k=1}^{n^2} s(k)
$$

where $s(k), k=1,2,...n^3$ is the decreasing rearrangement of $|x_ia_jb_k|,
i,j,k=1,2,...,n$. We apply Lemma 6 again with $s=n^2$ and the Orlicz
function N such that 

$$
N^{*}(\frac{1}{n^2} \sum_{j=1}^{l} t(j)) = \frac{l}{n^2}
\qquad l=1,2,...,n^2
$$

where $t(j), j=1,2,...,n^2$ is the decreasing rearrangement of
$|a_{i}b_{k}|, i,k=1,2,...,n$  and such that $N^{*}$ is an affine
function between the given values. We get

$$
\frac{1}{2} \frac{(n-1)^2}{n^2+(n-1)^2}\parallel x \parallel_{N} \leq
(\sum_{i=1}^{n}|x_i a_{\pi(i)}|^2)^{\frac {1}{2}} \leq
\frac{2}{c_n}\parallel x \parallel_{N} 
$$

We have for some integers $k_i$ with $k_i \leq n$ and 
$\sum k_i = ln$

$$
{N^{*}}^{-1}(\frac{l}{n}) = \frac{1}{n^2} \sum_{j=1}^{ln} t(j)
=\frac{1}{n^2} \sum_{i=1}^{n} a_{i} \sum_{j=1}^{k_{i}} b_{j}
=\frac{1}{n^2} \sum_{i=1}^{n} a_{i} \sqrt{nk_{i}}
$$

Since $a_1 \geq a_2 \geq \hdots \geq a_n \geq 0$ we also have
$k_1 \geq k_2 \geq \hdots \geq k_n$. Therefore we get

$$
\split
{N^{*}}^{-1}(\frac{l}{n}) &\leq \frac{1}{n^{\frac{3}{2}}}
 (\sqrt{k_1} \sum_{i=1}^{l} a_{i} +\sum_{i=l+1}^{n} a_{i} \sqrt{k_i}) \\
&\leq \frac{1}{n^{\frac{3}{2}}}
(|\sum_{i=1}^{l} a_{i}|^{2} + l \sum_{i=l+1}^{n} |a_{i}|^{2})^{\frac{1}{2}}
(k_{1}+ \frac{1}{l}\sum_{i=l+1}^{n} k_{i})^{\frac{1}{2}} \\
&\leq \frac{\sqrt{2}}{n}
(|\sum_{i=1}^{l} a_{i}|^{2} + l \sum_{i=l+1}^{n} |a_{i}|^{2})^{\frac{1}{2}} \\
&=\sqrt{2} {M^{*}}^{-1}(\frac{l}{n}) \\
\endsplit
$$

We get immediately that

$$
{N^{*}}^{-1}(\frac{l}{n}) = \frac{1}{n^2} \sum_{j=1}^{ln} t(j)
\geq \frac{1}{n} \sum_{i=1}^{l} a_{i}
$$

and as in [K--S 2]

$$
{N^{*}}^{-1}(\frac{l}{n}) \geq \frac{\sqrt{l}}{2n}
( \sum_{i=1}^{l} |a_{i}|^{2})^{\frac{1}{2}} 
$$ 

Therefore we have 

$$ 
 {M^{*}}^{-1}(\frac{l}{n})=
\frac{1}{n}(|\sum_{i=1}^{l} a_{i}|^{2} +
l \sum_{i=l+1}^{n} |a_{i}|^{2})^{\frac{1}{2}} 
\leq \sqrt{5} {N^{*}}^{-1}(\frac{l}{n})
$$

Altogether we have for $l=1,2,...,n$

$$ 
\frac{1}{\sqrt{2}}{N^{*}}^{-1}(\frac{l}{n}) \leq
 {M^{*}}^{-1}(\frac{l}{n}) \leq 
\sqrt{5} {N^{*}}^{-1}(\frac{l}{n})
$$

Since $M^{*}$ and $N^{*}$ are affine function for the other values the above
inequalities extend to arbitrary values and we get therefore

$$
\frac{1}{\sqrt{2}} \parallel x \parallel_{N}
\leq  \parallel x \parallel_{M}
\leq  \sqrt{5} \parallel x \parallel_{N}  
$$

\qed
\enddemo
\vskip 1cm

\proclaim{\smc Lemma 7}
Let H be a concave, increasing function on [0,1] that is twice continuously
differentiable on (0,1], continuous on [0,1] and satisfies $H(0)=0$. Assume
that $(\frac{H(t)}{t})' \ne 0$ for all $t \in (0,1]$. Then

$$
f(t)=-\frac{1}{2} \int_{t}^{1} \frac{H''(s)}{\sqrt{H(s)-sH'(s)}} ds
+\sqrt{H(1)} - \sqrt{H(1)-H'(1)}
\tag 5
$$

is a nonnegative, decreasing, differentiable function on (0,1] such that
$\int_{0}^{1} f(t)dt$ is finite and such that we have for all $t \in [0,1]$

$$
H(t)=(\int_{0}^{t} f(s)ds)^{2} + t \int_{t}^{1} |f(s)|^{2}ds
$$

\endproclaim
\vskip 1cm

\proclaim{\smc Lemma 8}
Let H be a concave, increasing function on [0,1] that is twice continuously
differentiable on (0,1], continuous on [0,1] and satisfies $H(0)=0$.
Moreover, assume that $(\frac{H(t)}{t})' \ne 0$ for all $t \in (0,1]$. Then
we have \par
\vskip 3mm

(i) $\lim_{t \to 0} t(-\frac{d}{dt}(\frac{H(t)}{t}))^{\frac{1}{2}}=0$ \par
\vskip 3mm

(ii) The function f given by (5) is well defined, nonnegative, decreasing,
and differentiable. \par
\vskip 3mm

(iii) $\lim_{t \to 0} tf(t)=0$ \par

\endproclaim
\vskip 1cm

\demo{Proof of Lemma 8}
(i) 

$$
\lim_{t \to 0} t(-\frac{d}{dt}(\frac{H(t)}{t}))^{\frac{1}{2}}= 
\lim_{t \to 0} t(\frac{H(t)}{t^2}-\frac{H'(t)}{t})^{\frac{1}{2}}=
\lim_{t \to 0}(H(t)-tH'(t))^{\frac{1}{2}}
$$

We use that $0 \leq H'(t) \leq \frac{H(t)}{t}$. \par
(ii) Because of $\frac{d}{dt}(\frac{H(t)}{t}) \ne 0$ and the concavity
of H we have $H(t)-tH'(t) > 0$. Again, by the concavity of H the integrand
is a nonpositive function and therefore f is nonnegative and decreasing.
\par
(iii) We have that

$$
\frac{1}{t} \frac{d}{dt} \left(t \sqrt{-\frac{d}{dt}(\frac{H(t)}{t})}\right)=
-\frac{H''(t)}{2\sqrt{H(t)-tH'(t)}}
\tag 6
$$

Integration by parts gives us

$$
tf(t)=t \left[\sqrt{-\frac{d}{ds}(\frac{H(s)}{s})}\right]_{t}^{1}
+t\int_{t}^{1} \frac{1}{s} \sqrt{-\frac{d}{ds}(\frac{H(s)}{s})}ds
+t(\sqrt{H(1)}-\sqrt{H(1)-H'(1)})
$$

The first summand tends to 0 because of (i) and the third trivially.
The second summand also tends to 0: If the integral is bounded this is
trivial. If the integral is not bounded we apply l'H\^opital's rule
and (i).
\qed
\enddemo
\vskip 1cm

\demo{Proof of Lemma 7}
In general f is unbounded in a neighborhood of 0.

$$
\int_{0}^{t} f(s)ds= \lim_{\epsilon \to 0}\{[sf(s)]_{\epsilon}^{t}
-\int_{\epsilon}^{t}sf'(s)ds\}
$$

By Lemma 8(ii), the definition (5) of f, and (6) we get

$$
\int_{0}^{t} f(s)ds=tf(t)+t\sqrt{-\frac{d}{dt}(\frac{H(t)}{t})}
-\lim_{\epsilon \to 0}
\epsilon \sqrt{-\frac{d}{d\epsilon}(\frac{H(\epsilon)}{\epsilon})}
$$

By Lemma 7(i) we get

$$
\int_{0}^{t} f(s)ds=tf(t)+t\sqrt{-\frac{d}{dt}(\frac{H(t)}{t})} 
\tag 7
$$

or

$$
\{\frac{1}{t} \int_{0}^{t} f(s)ds-f(t)\}^2=-\frac{d}{dt}(\frac{H(t)}{t})  
$$

Therefore we have

$$
\frac{H(x)}{x}-H(1)=-\int_{x}^{1} \frac{d}{dt}(\frac{H(t)}{t})dt
=\int_{x}^{1}|\frac{1}{t} \int_{0}^{t} f(s)ds - f(t)|^{2} dt
$$

With

$$
\frac{d}{dt}\{\frac{1}{t}(\int_{0}^{t}f(s)ds)^{2}+\int_{t}^{1}|f(s)|^{2}ds\}=
-|\frac{1}{t}\int_{0}^{t}f(s)ds-f(t)|^{2}
$$

we get

$$
\frac{H(x)}{x}-H(1)=\frac{1}{x}(\int_{0}^{x}f(s)ds)^{2}
+\int_{x}^{1}|f(s)|^{2}ds - (\int_{0}^{1}f(s)ds)^{2}
\tag 8
$$

By (7) we have 

$$
\int_{0}^{1}f(s)ds=f(1)+\sqrt{H(1)-H'(1)}
$$

By the definition (5) of f we get $f(1)=\sqrt{H(1)}-\sqrt{H(1)-H'(1)}$
and therefore

$$
\int_{0}^{1} f(s)ds= \sqrt{H(1)}
$$

Thus we obtain from(8)

$$
\frac{H(x)}{x}-H(1)=\frac{1}{x}(\int_{0}^{x}f(s)ds)^{2}
+\int_{x}^{1}|f(s)|^{2}ds -H(1)
$$

Or

$$
H(x)=(\int_{0}^{x} f(s)ds)^{2}+x \int_{x}^{1}|f(s)|^{2}ds
$$

\qed
\enddemo 
\vskip 1cm

\demo{Proof of Theorem 2}
Since M is strictly convex ${M'}^{-1}$ exists and ${M^{*}}'(t)={M'}^{-1}(t)$.
Since M is twice differentiable so is ${M^{*}}^{-1}$. Since $M(\sqrt t)$ is
strictly concave $({M^{*}}^{-1}(t))^2$ is also strictly concave. Therefore,

$$
0 > (({M^{*}}^{-1})^2)'(s)-\frac{({M^{*}}^{-1})^2(s)}{s}
= s \frac{d}{ds}(\frac{({M^{*}}^{-1}(s))^2}{s})
$$

We put $H(t)=({M^{*}}^{-1}(t))^2$ and apply
Lemma 7. Therefore $a_1, a_2,...., a_n$ given by (3) is a
positive, decreasing sequence with

$$
a_{l}=n \int_{\frac{l-1}{n}}^{\frac{l}{n}} f(s)ds
$$

We get

$$
\split
{M^{*}}^{-1}(\frac{l}{n})
&=((\int_{0}^{\frac{l}{n}} f(s)ds)^{2}+\frac{l}{n}
(\int_{\frac{l}{n}}^{1} |f(s)|^{2}ds))^{\frac{1}{2}} \\
&=((\frac{1}{n} \sum_{j=1}^{l} a_{j})^{2} + \frac{l}{n}
(\sum_{j=l}^{n-1} \int_{\frac{j}{n}}^{\frac{j+1}{n}} |f(s)|^{2}ds)
^{\frac{1}{2}}
\endsplit
$$

Since $f(\frac{l}{n}) \leq a_{l} \leq f(\frac{l-1}{n})$ we get

$$
{M^{*}}^{-1}(\frac{l}{n}) \leq
((\frac{1}{n} \sum_{j=1}^{l} a_{j})^{2} + \frac{l}{n}
(\frac{1}{n} \sum_{j=l}^{n-1} {a_{j}}^{2}))^{\frac{1}{2}}
$$

$${M^{*}}^{-1}(\frac{l}{n}) \geq
((\frac{1}{n} \sum_{j=1}^{l} a_{j})^{2} + \frac{l}{n}
(\frac{1}{n} \sum_{j=l+1}^{n} {a_{j}}^{2}))^{\frac{1}{2}}
$$

Now it is left to apply Theorem 1.
\qed
\enddemo

\enddocument
\bye